\documentclass[11pt]{article}
\usepackage{latexsym,a4wide}
\usepackage{amssymb}
\usepackage{amsmath}
\usepackage{color}
\usepackage{graphicx}
\usepackage{amsthm}
\usepackage{indentfirst}
\allowdisplaybreaks

\newcommand \bel {\begin{equation}\label}
\newcommand \ee {\end{equation}}
\newcommand \be {\begin{equation}}

\newcommand \RR {\mathbb R}

\newcommand \LL {\mathbb L}

\newcommand \del \partial

\newcommand \bei {\begin{itemize}}
\newcommand \eei {\end{itemize}}

\newtheorem{theorem}{\color{black}\indent Theorem}[section]

\newtheorem{proposition}{\color{black}\indent Proposition}[section]

\newtheorem{remark}{\color{black}\indent Remark}[section]

\begin{document}
\large
\title{Two family of explicit blowup solutions for $3$D incompressible Navier-Stokes equations }
\author{
{\sc Weiping Yan}\thanks{School of Mathematics, Xiamen University, Xiamen 361000, P.R. China. Email: yanwp@xmu.edu.cn.}
\thanks{Laboratoire Jacques-Louis Lions, Sorbonne Université, 4, Place Jussieu, 75252 Paris, France.}
}
\date{June 26, 2018}

\maketitle

\begin{abstract}  
In this paper, we give two family of explicit blowup solution for $3$D incompressible Navier-Stokes equations in $\RR^3$. Here one family of solutions admit the smooth initial data, and the initial data of another family of solutions are not smooth. The energy of those solutions is infinite.
\end{abstract}



\section{Introduction and main results} 
\setcounter{equation}{0}

The $3$D incompressible Navier-Stokes equations describes the motion of ideal incompressible fluid, which takes the following form
\bel{E1-1}
\aligned
&\textbf{v}_t+\textbf{v}\cdot\nabla \textbf{v}=-\nabla P+\nu\triangle \textbf{v},\\
&\nabla\cdot\textbf{v}=0,
\endaligned
\ee
where $\textbf{v}(t,x):[0,T^*)\times\RR^3\rightarrow\RR^3$ denotes the $3$D velocity field of the fluid, $P(t,x):[0,T^*)\times\RR^3\rightarrow\RR$ stands for the pressure in the fluid. $\nu>0$ denotes the kinematic viscosity.
The divergence free condition in second equations of (\ref{E1-1}) guarantees the incompressibility of the fluid.

We supplement the $3$D incompressible Navier-Stokes equations (\ref{E1-1}) with the initial data
$$
\textbf{v}(0,x)=\textbf{v}_0(x),\quad x\in\RR^3.
$$
It is easy to check that the functions (\ref{E1-5}) is a solution of $3$D solution of incompressible Navier-Stokes equations (\ref{E1-1}). 
The $3$D incompressible Navier-Stokes equations satisfies the energy identity
$$
\|\textbf{v}(t,x)\|^2_{\LL^2{(\RR^3})}+2\nu\int_0^t\|\nabla\textbf{v}(s,x)\|_{\LL^2(\RR^3)}=\|\textbf{v}(0,x)\|^2_{\LL^2{(\RR^3})}.
$$

We state the main result of this paper.

\begin{theorem}
\begin{itemize}
Let constant $T^*>0$ be the maximal existence time and the kinematic viscosity $\nu>0$.
The 3D incompressible Navier-Stokes equations (\ref{E1-1}) admits two family of explicit blowup solutions as follows

\item
One family of explicit solutions is
\bel{E1-7}
\textbf{v}(t,x)=\Big(v_1(t,x),v_2(t,x),v_3(t,x)\Big)^T,\quad (t,x)\in[0,T^*)\times\RR^3,
\ee
where
$$
\aligned
&v_1(t,x):={ax_1\over T^*-t}+{kx_2\over x_1^2+x_2^2},\\
&v_2(t,x):={ax_2\over T^*-t}-{kx_1\over x_1^2+x_2^2},\\
&v_3(t,x):=-{2ax_3\over T^*-t},
\endaligned
$$
with the initial data
\bel{E1-8}
\textbf{v}(0,x)=\Big({ax_1\over T^*}+{kx_2\over x_1^2+x_2^2},~{ax_2\over T^*}-{kx_1\over x_1^2+x_2^2},~-{2ax_3\over T^*}\Big)^T.
\ee

\item
Another family of explicit soutions is
\bel{E1-7R1}
\textbf{v}(t,x)=\Big(v_1(t,x),v_2(t,x),v_3(t,x)\Big)^T,\quad (t,x)\in[0,T^*)\times\RR^3,
\ee
where
$$
\aligned
&v_1(t,x):={ax_1\over T^*-t}+kx_2(T^*-t)^{2a},\\
&v_2(t,x):={ax_2\over T^*-t}-kx_1(T^*-t)^{2a},\\
&v_3(t,x):=-{2ax_3\over T^*-t},
\endaligned
$$
with the smooth initial data
\bel{E1-8R1}
\textbf{v}(0,x)=\Big({ax_1\over T^*}+kx_2(T^*)^{2a},~{ax_2\over T^*}-kx_1(T^*)^{2a},~-{2ax_3\over T^*}\Big)^T.
\ee
\end{itemize}
\end{theorem}

\begin{remark}
It follows from (\ref{E1-7}) and (\ref{E1-7R1}) that there is self-smilar singularity in $x_3$ direction, that is, ${2ax_3\over T^*-t}$ for $a\in\RR/\{0\}$.
\end{remark}

\begin{remark}

On one hand, it follows from (\ref{E1-7}) that
$$
\aligned
&\nabla v_1(t,x)=\Big({a\over T^*-t}-{2kx_1x_2\over(x_1^2+x_2^2)^2},~{k(x_1^2-x_2^2)\over (x_1^2+x_2^2)^2},~0\Big)^T,\\
&\nabla v_2(t,x)=\Big({k(x_1^2-x_2^2)\over (x_1^2+x_2^2)^2},~{a\over T^*-t}+{2kx_1x_2\over(x_1^2+x_2^2)^2},~0\Big)^T,\\
&\nabla v_3(t,x)=\Big(0,~0,~-{2a\over T^*-t}\Big)^T,
\endaligned
$$
which means that
$$
div(v_i)|_{x=x_0}=\infty,\quad as\quad t\rightarrow (T^*)^-,
$$
for a fixed point $x_0\in\RR^3$. Here one can see the initial data is not smooth from (\ref{E1-8}).

On the other hand, it follows from (\ref{E1-7R1}) that
$$
\aligned
&\nabla v_1(t,x)=\Big({a\over T^*-t},~k(T^*)^{2a},~0\Big)^T,\\
&\nabla v_2(t,x)=\Big(-k(T^*)^{2a},~{a\over T^*-t},~0\Big)^T,\\
&\nabla v_3(t,x)=\Big(0,~0,~-{2a\over T^*-t}\Big)^T,
\endaligned
$$
which means that
$$
div(v_i)|_{x=x_0}=\infty,\quad as\quad t\rightarrow (T^*)^-,
$$
for a fixed point $x_0\in\RR^3$. Here one can see the initial data is smooth from (\ref{E1-8R1}). But the initial data goes to infinity as $x\rightarrow\infty$.

Hence those two family of blowup solutions (\ref{E1-7}) and (\ref{E1-7R1}) are not finite energy solutions.
\end{remark}

\begin{remark}
It is well-known that the equation for pressure $P$ is
$$
-\triangle P=\sum_{i,j=1}^3{\del v_i\over\del x_j}{\del v_j\over\del x_i}.
$$
Hence by (\ref{E1-6}), we can get the pressure $P$.
\end{remark}

Let $\textbf{e}_r$, $\textbf{e}_{\theta}$ and $\textbf{e}_z$ be the cylindrical coordinate
system,
\bel{E1-0}
\aligned
&\textbf{e}_r=({x_1\over r},{x_2\over r},0)^T,\\
&\textbf{e}_{\theta}=({x_2\over r},-{x_1\over r},0)^T,\\
&\textbf{e}_z=(0,0,1)^T,
\endaligned
\ee
where $r=\sqrt{x_1^2+x_2^2}$ and $z=x_3$.

Then we can rewrite explicit blowup solutions (\ref{E1-7}) and (\ref{E1-7R1}) into two axisymmetric blowup solutions.

\begin{theorem}
\begin{itemize}
Let $T^*>0$ be a constant,  the kinematic viscosity $\nu>0$, and $\textbf{e}_r,\textbf{e}_{\theta},\textbf{e}_z$ are defined in (\ref{E1-0}), $r=\sqrt{x_1^2+x_2^2}$ and $z=x_3$.
The 3D incompressible Navier-Stokes equations (\ref{E1-1}) admits two family of explicit blowup axisymmetric solutions: 

\item One family of explicit blowup axisymmetric solutions is
\bel{E1-6}
\textbf{v}(t,x)=v^r(t,r,z)\textbf{e}_r+v^{\theta}(t,r,z)\textbf{e}_{\theta}+v^z(t,r,z)\textbf{e}_z,\quad (t,x)\in[0,T^*)\times\RR^3,
\ee
where
$$
\aligned
&v^r(t,r,z)={ar\over T^*-t},\\
&v^{\theta}(t,r,z)={k\over r},\\
&v^z(t,r,z)=-{2az\over T^*-t},
\endaligned
$$
where constants $a,k\in\RR/\{0\}$.

\item Another family of explicit blowup axisymmetric solutions is
\bel{E1-6r1}
\textbf{v}(t,x)=v^r(t,r,z)\textbf{e}_r+v^{\theta}(t,r,z)\textbf{e}_{\theta}+v^z(t,r,z)\textbf{e}_z,\quad (t,x)\in[0,T^*)\times\RR^3,
\ee
where
$$
\aligned
&v^r(t,r,z)={ar\over T^*-t},\\
&v^{\theta}(t,r,z)= kr(T^*-t)^{2a},\\
&v^z(t,r,z)=-{2az\over T^*-t},
\endaligned
$$
where constants $a,k\in\RR/\{0\}$.
\end{itemize}
\end{theorem}

\begin{remark}
Comparing with two solutions (\ref{E1-6}) and (\ref{E1-6r1}), there is only different component in $\textbf{e}_{\theta}$ direction. But this causes the smoothness of initial data.
\end{remark}

In 1934, Leray \cite{L} showed that the $3$D incompressible Navier-Stokes equations admit global-forward-in-time weak solution of the initial value problem. It is not known, however, whether the Leray solutions are unique, or whether the solutions for suitable general smooth initial data remain smooth or blowup \cite{F}. 
After that, there are many papers concerns with the well-posedness of weak solutions or blowup of solution for this problem. One can see \cite{B,Ca,CP,CKN,HL,JS,NRS,S,T} for more details. 
For the existence of explicit solutions of incompressible fluids, let us note that the first nontrivial example of hidden symmetries connected with reduction of PDEs was found by Kapitanskiy \cite{K,K0} for the Navier-Stokes equations.
Constantin \cite{Con} found a class of smooth, mean zero initial data for which the solution of $3$D Euler equations becomes infinite in finite time, meanwhile, he gave an explicit formulas of solutions for the $3$D Euler equations by reducing this equations into a local conservative Riccati system in two-dimensional basic square. 
Due to space limitations, here we can not list all of interesting results.
In this paper, we study the explicit blowup solutions for the incompressible Navier-Stokes equations (\ref{E1-1}).

\section{Proof of Thereom 1.1}

We now derive the $3$D incompressible Navier-Stokes equations with axisymmetric velocity field in the cylindrical coordinate.
The $3$D velocity field $\textbf{v}(t,x)$ is called axisymmetric if it can be written as
$$
\textbf{v}(t,x)=v^r(t,r,z)\textbf{e}_r+v^{\theta}(t,r,z)\textbf{e}_{\theta}+v^z(t,r,z)\textbf{e}_z,
$$
where $v^r$, $v^{\theta}$ and $v^z$ do not depend on the $\theta$ coordinate.

Thus in the cylindrical coordinates, the $3$D Navier-Stokes equations (\ref{E1-1}) with axisymmetric velocity field can be reduced into a system as follows
\bel{E1-2}
\aligned
&v_t^{\theta}+v^{r}v_r^{\theta}+v^zv_z^{\theta}=\nu\Big(\triangle-{1\over r^2}\Big)v^{\theta}-{v^rv^{\theta}\over r},\\
&\omega_t^{\theta}+v^r\omega_r^{\theta}+v^z\omega_z^{\theta}=\nu\Big(\triangle-{1\over r^2}\Big)\omega^{\theta}+{2\over r}v^{\theta}v_z^{\theta}+{1\over r}v^r\omega^{\theta},\\
&-\Big(\triangle-{1\over r^2}\Big)\phi^{\theta}=\omega^{\theta},
\endaligned
\ee
where the radial and angular velocity fields $v^r$ and $v^{\theta}$ are recovered from $\phi^{\theta}$ based on the Biot-Savart law
\bel{E1-3}
v^r=-\del_z\phi^{\theta},\quad v^z={1\over r}\del_r(r\phi^{\theta}).
\ee

The incompressibility condition becomes
\bel{E1-4}
\del_r(rv^r)+\del_z(rv^z)=0.
\ee

The following result gives two family of explict self-similar blowup solutions for system (\ref{E1-2})-(\ref{E1-3}) with the  incompressibility condition (\ref{E1-4}).

\begin{proposition}
Let $T^*>0$ be a constant and $\nu>0$. System (\ref{E1-2})-(\ref{E1-3}) with the incompressibility condition (\ref{E1-4}) admits two family of explicit solutions $(v^r,v^{\theta},v^z,\phi^{\theta},\omega^{\theta})$, which have the form
\bel{E1-5}
\aligned
&v^r(t,r,z)={ar\over T^*-t},\\
&v^{\theta}(t,r,z)={k\over r},\quad or\quad kr(T^*-t)^{2a},\\
&v^z(t,r,z)=-{2az\over T^*-t},\\
&\phi^{\theta}=-{arz\over T^*-t},\\
&\omega^{\theta}=0,
\endaligned
\ee
where constants $a,k\in\RR/\{0\}$.
\end{proposition}

Since the angular components $v^{\theta}(t,r,z)$, $\omega^{\theta}(t,r,z)$ and $\phi^{\theta}(t,r,z)$ can be seen as odd functions of $r$ \cite{LW},  we follow \cite{HL} to introduce the following transformation
\bel{E1-4R1}
v_1={v^{\theta}\over r},\quad \omega_1={\omega^{\theta}\over r},\quad \phi_1={\phi^{\theta}\over r},
\ee
then system (\ref{E1-2}) is reduced into 
\bel{E1-2R1}
\del_tv_1+v^{r}\del_rv_1+v^z\del_zv_1=\nu\Big(\del_r^2+{3\over r}\del_r+\del_z^2\Big)v_1+2v_1\del_z\phi_1,
\ee
\bel{E1-2R2}
\del_t\omega_1+v^r\del_r\omega_1+v^z\del_z\omega_1=\nu\Big(\del_r^2+{3\over r}\del_r+\del_z^2\Big)\omega_1+\del_z(v_1^2),
\ee
\bel{E1-2R3}
-\Big(\del_r^2+{3\over r}\del_r+\del_z^2\Big)\phi_1=\omega_1,
\ee
with the Biot-Savart law given by
\bel{E1-3R1}
v^r=-r\del_z\phi_1,\quad v^z=2\phi_1+r\del_r\phi_1.
\ee
The incompressibility condition (\ref{E1-4}) is invariant . 

\paragraph{Proof of Proposition 2.1.}
The idea of finding explicit blowup solutions for system (\ref{E1-2R1})-(\ref{E1-3R1}) with the incompressibility condition (\ref{E1-4}) comes from \cite{Yan1,Yan2,Yan0}. This is based on the observation on the structure of system (1.3)-(1.4) and incompressibility condition (1.5). 

Let $a\neq0$ be a unknown constant in $\RR$. 
By the the structure of incompressibility condition (\ref{E1-4}), we assume that
\bel{E2-1}
\aligned
&v^r(t,r,z)={ar\over T^*-t},\\
&v^z(t,r,z)=-{2az\over T^*-t}.
\endaligned
\ee
It is easy to see $v^r(t,r,z)$ and $v^z(t,r,z)$ given in (\ref{E2-1}) satisfies the incompressibility condition (\ref{E1-4}).

Substituting (\ref{E2-1}) into (\ref{E1-3R1}), we get
$$
\aligned
&-r\del_z\phi_1={ar\over T^*-t},\\
&2\phi_1+r\del_r\phi_1=-{2az\over T^*-t},
\endaligned
$$
which gives that 
\bel{E2-2}
\phi_1(t,r,z)=-{az\over T^*-t}.
\ee

Substituting (\ref{E2-2}) into (\ref{E1-2R3}), there is
\bel{E2-3}
\omega_1=0.
\ee

By (\ref{E2-3}), it follows from (\ref{E1-2R2}) that
$$
\del_z(v_1)^2=0,
$$
from which, we assume that 
\bel{E2-4}
v_1(t,r,z)={kr^p\over (T^*-t)^{\beta}},
\ee
where $k\neq0$, $p$ and $\beta$ are three unknown constants in $\RR$.

We substitute (\ref{E2-1})-(\ref{E2-4}) into equation (\ref{E1-2R1}), we have
$$
{k\beta r^p\over (T^*-t)^{\beta+1}}+{akpr^p\over (T^*-t)^{\beta+1}}={\nu\Big(kp(p-1)r^{p-2}+3kpr^{p-2}\Big)\over (T^*-t)^{\beta}}-{2akr^p\over (T^*-t)^{\beta+1}}.
$$

If above equation holds, then we should require 
$$
\aligned
&kp(p-1)+3kp=0,\\
&k\beta+akp+2ak=0,
\endaligned
$$
which gives that
$$
p=-2, \quad \beta=0,
$$
or
$$
p=0, \quad \beta=-2a,
$$
Thus by (\ref{E2-4}), we get
\bel{E2-5}
v_1(t,r,z)={k\over r^2},
\ee
or 
\bel{E2-6}
v_1(t,r,z)=k(T^*-t)^{2a}.
\ee

Hence it follows from (\ref{E1-4R1}) that
$$
\aligned
&v^r(t,r,z)={ar\over T^*-t},\\
&v^{\theta}(t,r,z)={k\over r},\quad or\quad kr (T^*-t)^{2a},\\
&v^z(t,r,z)={-2az\over T^*-t},\\
&\phi^{\theta}(t,r,z)=-{arz\over T^*-t},\\
&\omega^{\theta}=0,
\endaligned
$$
which are the solution of system (\ref{E1-2})-(\ref{E1-3}) with the incompressibility condition (\ref{E1-4}). Here constants $a,k\in\RR/\{0\}$.

\paragraph{Proof of Theorem 1.2.}
Since
$$
\textbf{v}(t,x)=v^r(t,r,z)\textbf{e}_r+v^{\theta}(t,r,z)\textbf{e}_{\theta}+v^z(t,r,z)\textbf{e}_z,\quad (t,x)\in[0,T^*)\times\RR^3,
$$
we can obtain a family of explicit blowup axisymmetric solutions for $3$D incompressible Navier-Stokes equations by noticing that $\textbf{e}_r,\textbf{e}_{\theta},\textbf{e}_z$ are defined in (\ref{E1-0}), $r=\sqrt{x_1^2+x_2^2}$ and $z=x_3$, and
$$
\aligned
&v^r(t,r,z)={ar\over T^*-t},\\
&v^{\theta}(t,r,z)={k\over r},\quad or\quad kr(T^*-t)^{2a},\\
&v^z(t,r,z)=-{2az\over T^*-t},
\endaligned
$$
where constants $a,k\in\RR/\{0\}$. 

Direct computations give the vorticity vector: 
$$
\omega(t,r,z)=\textbf{0},\quad or\quad kr(T^*-t)^{2a}\textbf{e}_z,
$$
where constant $k\in\RR/\{0\}$.

\paragraph{Proof of Theorem 1.1.}
By directly computations, we can obtain two family of explicit blowup solutions from (\ref{E1-6}) and (\ref{E1-6r1}).


\

\

\

\textbf{Acknowledgments.} 

The author expresses his sincerely thanks to the BICMR of Peking University and Professor Gang Tian for constant support and encouragement.
The author also expresses his sincerely thanks to Prof. V. Sver\'{a}k for informed the paper \cite{Con,K,K0} and his suggestions. 
The author is supported by NSFC No 11771359.

\end{document}